\documentclass[12pt]{article}
\usepackage[cp1251]{inputenc}
\usepackage[dvips]{graphicx}
\usepackage{amssymb,amsmath}
\usepackage{amsthm}
\usepackage{epsf}
\author{Alexander I. Efimov\thanks{e-mail: efimov13@yandex.ru},\\ \emph{Faculty of Mechanics
and Mathematics of Moscow State University;}\\ \emph{Independent University of Moscow}}
\title{The asymptotics for the number of real
roots of the Bernoulli polynomials}
\date{}

\begin{document}

\maketitle

\smallskip
{\bf Abstract.} {\small A new short clear proof of the asymptotics
for the number $c_n$ of real roots of the Bernoulli polynomials
$B_n(x)$, as well as for the maximal root $y_n$:
$$ y_n=\frac{n}{2\pi e}+\frac{\ln(n)}{4\pi e}+O(1)\quad\text{and}\quad
c_n=\frac{2n}{\pi e}+\frac{\ln(n)}{\pi e}+O(1). $$}

\smallskip
{\bf I. Introduction.} In this paper we study the
behavior of real roots of the Bernoulli polynomials. For
motivations see [1]. The Bernoulli polynomials $ B_n(x) $ are
defined by the following identity in the ring $ \mathbb Q[x][[t]]
$.
$$
\frac{te^{tx}}{e^t-1}=\sum\limits_{n=0}^\infty
{\frac{B_n(x)}{n!}}t^n\text.
$$
Thus, $ B_0(x)=1 $, $ B_1(x)=x-\frac12 $, $ B_2(x)=x^2-x+\frac16$,...
An equivalent definition is the following: $B_0(x)=1$, for each $n$
the polynomial $ B_{n+1}(x) $ is a primitive of the polynomial
$(n+1)B_n(x) $ such that
$$ \int\limits_0^1B_{n+1}(x)dx=0.
$$

We are interested in the asymptotics for the maximal real root $y_n$ and
for the number $c_n$ of real roots of polynomials $B_n(x)$ as $n\to\infty$. In
this paper the following results are proved.


\smallskip
{\bf Theorem 1.} {\it As $n\to\infty$ }
$$ y_n=\frac{n}{2\pi e}+\frac{\ln(n)}{4\pi e}+O(1)\text. $$
By $]x[$ we denote the minimal integer greater or equal to $x$.

\smallskip
{\bf Theorem 2.} {\it The following relation is fulfilled }
$$ c_{4k+1}=4]y_{4k+1}[-3. $$

As a corollary of Theorems 1 and 2 we obtain an asymptotics for
the number of real roots of the Bernoulli polynomials $ B_n(x) $
as $n\to\infty$.

\smallskip
{\bf Theorem 3.} {\it As $n\to\infty$}
$$ c_n=\frac{2n}{\pi e}+\frac{\ln(n)}{\pi e}+O(1). $$

Theorems 1 and 3 are stronger than K. Inkeri's results [3]:
$$y_n=\frac{n}{2\pi e}+o(n)\quad,\quad c_n=\frac{2n}{\pi e}+o(n).$$

H. Delange proved Theorems 1 and
3 in papers [4],[5] and showed
estimates for $O(1)$. Theorem 2 is actually a corollary of Inkeri's results.
The proof given in this paper is simpler and conceptually clearer.
Our proofs are in the following steps (see the details below).

Throughout all proofs we use the identities [1, \S 15]

(1)\qquad $B_n(1-x)={(-1)^n}B_n(x)$;

(2)\qquad $2^{1-n}B_n(2x)=B_n\left(x\right)+B_n\left(x+\frac12\right)$;

(3)\qquad $B_n(x+1)-B_n(x)=nx^{n-1}$;

(4)\qquad $B_n'(x)=nB_{n-1}(x)$.

The proof of Theorem 1 is as follows.
First using these identities we prove (the Lemma below) that

for $n$ even $B_n(x)$ increases and decreases on $[0,1]$ as it is shown in figures 1, 3.

for $n$ odd $B_n(x)$ is positive or negative on $[0,1]$ as it is shown in figures 2, 4.

\smallskip
\begin{tabular}{c}
\epsfbox{a.3}\\
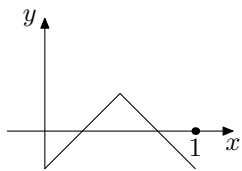
Figure 1.\\
$n\equiv 0\pmod 4$
\end{tabular}
\hfil
\begin{tabular}{c}
\epsfbox{a.4}\\
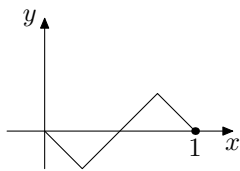
Figure 2.\\
$n\equiv 1\pmod 4$
\end{tabular}
\hfil
\begin{tabular}{c}
\epsfbox{a.1}\\
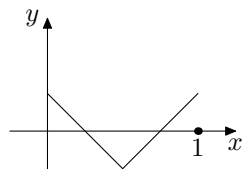
Figure 3.\\
$n\equiv 2\pmod 4$
\end{tabular}
\hfil
\begin{tabular}{c}
\epsfbox{a.2}\\
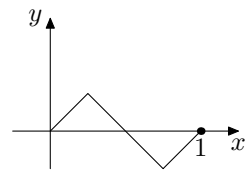
Figure 4.\\
$n\equiv 3\pmod 4$
\end{tabular}




\smallskip
Second, denote $d_n=]y_n[$.
Using the identities (3) and (4) together with Lemma we obtain

\smallskip
{\bf Statement 1.} $d_{n+1}\leq d_n+1$.

\smallskip
So $d_{4k}+i\geq d_{4k+i}\geq d_{4k+4}-4+i$ for $i\in\{1;2;3\}$ and
it suffices to find the asymptotics for $d_{4k}$ but not for $d_n$.

Using the identity (3) and integral estimates for
$\sum\limits_{i=1}^{m}i^{4k-1}$ we obtain

\smallskip
{\bf Statement 2.} $\sqrt[4k]{1-B_{4k}(0)}<d_{4k}<2+\sqrt[4k]{-B_{4k}(0)}.$

\smallskip
So it suffices to find the asymptotics for the radical $\sqrt[4k]{-B_{4k}(0)}$.

\smallskip
{\bf Statement 3.} $\sqrt[4k]{-B_{4k}(0)}=\dfrac{2k}{\pi e}+\dfrac{\ln(4k)}{4\pi
e}+O(1).$

\smallskip
We find the required asymptotics using the Euler formula [2, p.145] and the Stirling formula
$$\zeta(4k)=\frac{-B_{4k}(0)2^{4k-1}\pi^{4k}}{(4k)!}\quad\text{and}\quad
(4k)!\sim\left(\frac{4k}{e}\right)^{4k}\sqrt{8\pi k}.$$
Here $\zeta(4k)$ is the Riemann zeta-function of which we only use that
$1<\zeta(4k)<2$.
Thus Theorem 1 is proved (modulo the details).

Using Lemma and the identity (3) we obtain that the polynomial $B_{4k+1}(x)$ has exactly 2 roots
on the interval $[m;m+1)$ for $1\leq m\leq d_{4k+1}-1$ (taking multiplicity into account).
Then using (1) we prove Theorem 2.

Theorem 3 is proved using Theorems 1 and 2 and the inequality $$c_{n+1}\leq c_n+1.$$

The author is grateful for A.L.Gorodentsev, A.A.Karatsuba,
A.B.Skopenkov and A.V.Ustinov for
helpful and stimulating discussions.

\bigskip
\newpage {\noindent {\bf II. The details of the proofs.}}

\smallskip
{\bf Lemma.} {\it If $n\geq2$, then the following statements hold:

1. If $n\equiv2\text{(mod 4)}$, then $B_n(0)=B_n(1)>0>B_n(\frac12)$; moreover, the function $B_n(x)$ is
strictly decreasing on the segment $\left[0;\frac12\right]$,
strictly increasing on the segment $\left[\frac12;1\right]$.

2. If $n\equiv0\text{(mod 4)}$, then $B_n(0)=B_n(1)<0<B_n(\frac12)$; moreover, function $B_n(x)$ is
strictly increasing on segment $\left[0;\frac12\right]$, is
strictly decreasing on segment $\left[\frac12;1\right]$}.

3. If $n\equiv1\text{(mod 4)}$, then $B_n(0)=B_n(\frac12)=B_n(1)=0$; moreover, $B_n(x)<0$ for
$x\in\left(0;\frac12\right)$, $B_n(x)>0$ for $x\in\left(\frac12;1\right)$.

4. If $n\equiv3\text{(mod 4)}$, then $B_n(0)=B_n(\frac12)=B_n(1)=0$; moreover, $B_n(x)>0$ for
$x\in\left(0;\frac12\right)$, $B_n(x)<0$ for $x\in\left(\frac12;1\right)$.

\smallskip
{\bf Proof.} The proof is by induction over $n$.

For $n=2$ there is nothing to prove.

Suppose that the inductive hypothesis holds for each $m$ such that $2\leq m\leq n$.
Let us prove the Statement for $n+1$.

{\it The case $n\equiv2\text{(mod 4)}$.} Using (1) we get $B_{n+1}(0)=(-1)^{n+1}B_{n+1}(1)=-B_{n+1}(1)$. Using
(3) we get
$B_{n+1}(1)=B_{n+1}(0)+(n+1)0^n=B_{n+1}(0)$. Hence, $B_{n+1}(1)=B_{n+1}(0)=0$.
Further, $B_{n+1}(\frac12)=(-1)^{n+1}B_{n+1}(\frac12)=-B_{n+1}(\frac12)$,
$B_{n+1}(\frac12)=0$. Using (4) and the
inductive hypothesis, for $x\in\left]0;\frac12\right[$ we have
$B_{n+1}''(x)=(n+1)nB_n(x)<0$. So function $B_{n+1}(x)$ is convex
up on segment $\left[0;\frac12\right]$. Taking into account
$B_{n+1}(0)=B_{n+1}(\frac12)=0$, we get $B_{n+1}(x)>0$ for $x\in\left]0;\frac12\right[$.
For $x\in\left]\frac12;1\right[$
we have $B_{n+1}(x)=-B_n(1-x)<0$. The inductive step is proved in this case.

{\it The case $n\equiv0\text{(mod 4)}$.} Analogously to the previous case.

{\it The case $n\equiv1\text{(mod 4)}$.}
Using (1) we get $B_{n+1}(0)=(-1)^{n+1}B_{n+1}(1)=B_{n+1}(1)$.
Using (2) for $x=0$ we get $B_{n+1}(0)=2^n\left(B_{n+1}(0)+B_{n+1}\left(\frac12\right)\right)$.
Thus $B_{n+1}(\frac12)=\frac{1-2^n}{2^n}B_{n+1}(0)$.
So, the signs of $B_{n+1}(0)$ and $B_{n+1}(\frac12)$ are different.
Using (4) and the inductive hypothesis we obtain
$$B_{n+1}'(x)=(n+1)B_n(x)<0\quad{for}\quad x\in\left(0;\frac12\right)$$ and
$$B_{n+1}'(x)=(n+1)B_n(x)>0\quad{for}\quad x\in\left(\frac12;1\right).$$
So the function $B_{n+1}(x)$ is strictly decreasing on segment
$\left[0;\frac12\right]$, is strictly increasing on segment $\left[\frac12;1\right]$, and
$B_{n+1}(0)=B_{n+1}(1)>0>B_{n+1}(\frac12)$.
The inductive step is proved in this case.

{\it The case $n\equiv3\text{(mod 4)}$.} Analogously to the previous case.

The inductive step is completely proved. The Lemma is proved.\qed

\smallskip
{\bf Corollary of (3) and (4).} {\it The degree of the polynomial $B_n(x)$ is $n$,
the top coefficient is 1; for each $ x>0 $ the sequence $
\{B_n(x+k)\}_{k=0}^{\infty} $ is strictly increasing}.

\smallskip
{\bf Proof of Statement 1.}
If $d_{n+1}\leq 2$, then $d_{n+1}\leq d_n+1$, because $d_n\geq 1$.
Now assume that $d_{n+1}\geq 3.$

1) {\it The case $n\equiv0\text{, 1, 2(mod 4)}$.} For these values of $n$ we
will get that $d_{n+1}\leq d_n$. Using the Lemma we obtain
$B_{n+1}(1)\geq 0$.
So, using the Corollary we obtain $B_{n+1}(d_{n+1}-1)>0 $.
Thus $ B_{n+1}(d_{n+1})>0$.
So the values of the polynomial $B_{n+1}(x)$ at the
endpoints of the segment $[d_{n+1}-1;d_n]$ are positive.
We also have that $y_{n+1}\in [d_{n+1}-1;d_n]$ is the root of the polynomial.
So the polynomial $B_{n+1}(x)$ has at least 2 roots on the interval $]d_{n+1}-1;d_n[$
(taking multiplicity into account).
Then there is a root of the polynomial $B_{n+1}'(x)$ between these roots.
From (4) this root is a root of the polynomial $B_n(x)$.
Thus $y_n>d_{n+1}-1$.
Hence $d_{n+1}\leq ]y_n[=d_n$.

2) $n\equiv3\text{(mod 4)}$. From the Lemma we get
$B_{n+1}(0)<B_{n+1}(x)$ for each $x\in ]0;1[$.
So using (3) we obtain that
$B_{n+1}(m)<B_{n+1}(m+x)$ for $m\in \mathbb{N}$, $x\in ]0;1[$.
Hence, $B_{n+1}(d_{n+1}-1)<B_{n+1}(y_{n+1})=0$. So from the
Corollary we get $B_{n+1}(d_{n+1}-2)<0$. Moreover,
$B_{n+1}\left(d_{n+1}-\frac32\right)>0$, because
$B_{n+1}(\frac12)>0$.

Further, as in the previous case, on the interval
$\left[d_{n+1}-2;d_{n+1}-1\right]$ there are at least 2 roots of
the polynomial $B_{n+1}(x)$.
Now the conclusion proof is analogous to the previous case.

Thus, the Statement 1 is proved.\qed

\smallskip
{\bf Proof of Statement 2.}
{\it Lower estimation.} Since the top coefficient of the polynomial $B_{4k}$
is equal to 1, and each root does not exceed $d_{4k}$, it follows that $B_{4k}(d_{4k})\geq 0$.
Consequently, from (3) we get
$$-B_{4k}(0)\leq 4k\sum\limits_{i=1}^{d_{4k}-1}i^{4k-1}<4k\int\limits_1^{d_{4k}}x^{4k-1}dx=
d_{4k}^{4k}-1.$$

{\it Upper estimation.} As it is shown above, $B_{4k}(d_{4k}-1)<0$.
Therefore from (3) we get that
$$-B_{4k}(0)>4k\sum\limits_{i=1}^{d_{4k}-2}i^{4k-1}>4k\int\limits_{0}^{d_{4k}-2}x^{4k-1}dx=(d_{4k}-2)^{4k}$$

Thus $$
\sqrt[4k]{1-B_{4k}(0)}<d_{4k}<2+\sqrt[4k]{-B_{4k}(0)}\quad\qed
$$

\smallskip
{\bf Proof of Statement 3.}
By the Euler formula [2, p.145], $$\zeta(4k)=\frac{-B_{4k}(0)2^{4k-1}\pi^{4k}}{(4k)!}.$$
By the Stirling formula, as $k\to\infty$
$$(4k)!\sim\left(\frac{4k}{e}\right)^{4k}\sqrt{8\pi k}.$$
 Now we have:
$$\sqrt[4k]{-B_{4k}(0)}=\sqrt[4k]{\frac{\zeta(4k)(4k)!}{2^{4k-1}\pi^{4k}}}=
\left(\sqrt[4k]{2\zeta(4k)}-1\right)\frac1{2\pi}\sqrt[4k]{(4k)!}+\frac1{2\pi}
\sqrt[4k]{(4k)!}.$$
In the last sum the first summand is $O(1)$, because
$\zeta(s)<2$ and so $\sqrt[4k]{2\zeta(4k)}-1=O(1/k)$.
Further,
\begin{multline*}
\frac1{2\pi}\sqrt[4k]{(4k)!}+O(1)=\frac{2k}{\pi e}+\frac{2k}{\pi
e}\left(\sqrt[8k]{(2\pi+o(1))4k}-1\right)+O(1)=\\=\frac{2k}{\pi
e}+\frac{2k}{\pi e}\left(\sqrt[8k]{4k}-1\right)+O(1)=
\frac{2k}{\pi e}+\frac{2k}{\pi
e}\left(e^{\frac{\ln(4k)}{8k}}-1\right)+O(1)=\\=\frac{2k}{\pi
e}+\frac{2k}{\pi
e}\left(\frac{\ln(4k)}{8k}+O\left(\left(\frac{\ln(4k)}{8k}\right)^2\right)\right)+O(1)=
\frac{2k}{\pi e}+\frac{\ln(4k)}{4\pi e}+O(1).
\end{multline*}

Thus, $$
\sqrt[4k]{-B_{4k}(0)}=\frac{2k}{\pi e}+\frac{\ln(4k)}{4\pi e}+O(1).\quad\qed
$$


\smallskip
{\bf Proof of Theorem 2.}
From (1) we obtain that the set of roots of polynomial $B_n(x)$ is
symmetric with respect to the point $\frac12$.
So it suffices to consider only roots greater than $\frac12$.
From (4) and the Lemma we obtain $B_{4k+1}(1)=0$ and
$B_{4k+1}'(1)=(4k+1)B_{4k}(1)<0$.
Hence the polynomial $B_{4k+1}(x)$ has a root greater than 1, so
$d_{4k+1}\geq 2$.

Let $m$ be an integer such that $1\leq m\leq d_{4k+1}-1$.
Then from the Corollary we obtain $B_{4k+1}(m+1)>0$.
Hence $B_{4k+1}(m)>0$ for $m>1$. For $x\in[m+\frac12;m+1]$, from
the Lemma we obtain $B_{4k+1}(x-m)\geq 0$.
Therefore $B_{4k+1}(x)>0$.
So there are no roots of the polynomial $B_{4k+1}(x)$ on segment
$[m+\frac12;m+1]$. In particular,
$y_{4k+1}\in]d_{4k+1}-1;d_{4k+1}-\frac12[$. Thus
$B_{4k+1}(y_{4k+1}-d_{4k+1}+1+m)\leq 0$, because $B_{4k+1}(y_{4k+1})=0$.

So for each $1\leq m\leq d_{4k+1}-1$ we have that:

1) At some point on segment $]m;m+\frac12[$ the value of the function
$B_{4k+1}(x)$ is not positive.

2) $B_{4k+1}(m)\geq 0$.

3) $B_{4k+1}(m+\frac12)>0$ (because from the Lemma we obtain
$B_{4k+1}(\frac12)>0$).

4) $B_{4k+1}''(x)=(4k+1)4kB_{4k-1}(x)>0$ for $x\in]m;m+\frac12]$
(because from the Lemma we obtain $B_{4k-1}(x-m)>0 $).

So the polynomial $ B_{4k+1}(x) $ has exactly 2 roots
on segment $[m;m+\frac12]$ (taking multiplicity into account).
Thus
the polynomial $B_{4k+1}(x)$ has exactly $2(d_{4k+1}-1)$ roots greater than
$\frac12$.
It also has the same number of roots less than $\frac12$, and $\frac12$ is a
root of multiplicity one. So
$$c_{4k+1}=4d_{4k+1}-3.\quad\qed$$

\smallskip
{\bf Proof of Theorem 3.} From Theorems 1 and 2 we obtain
$$c_{4k+1}=\frac{2(4k+1)}{\pi e}+\frac{\ln(4k+1)}{\pi e}+O(1).$$
From (4) we get $B_{n+1}'(x)=(n+1)B_{n}(x)$. Therefore $$c_{n+1}\leq
c_n+1\text.
$$ Hence $c_{4k+1}+i\geq c_{4k+1+i}\geq
c_{4k+5}-4+i$ for $i\in\{1;2;3\}$. Thus
$$c_n=\frac{2n}{\pi e}+\frac{\ln(n)}{\pi e}+O(1).$$
Theorem 3 is proved.\qed

\bigskip
{\large References:}\\
\\
$ [1] $ V.V. Prasolov, Polynomials. M.:MCCME, 2001.\\
$ [2] $ J.-P. Serre, Arithmetic course. M.:Mir, 1972.\\
$ [3] $ Inkeri, K. The real roots of Bernoulli polynomials. Ann.
Univ. Turku. Ser. A I 37 (1959) 3-19.\\
$ [4] $ H. Delange, Sur les zeros reels des polynomes de
Bernoulli. C. R. Acad. Sc. Paris, s

erie I, 303 (1986), 539-542.\\
$ [5] $ H. Delange, Sur les zeros reels des polynomes de
Bernoulli. Ann. Inst. Fourier, Grenoble, 41, 2 (1991), 267-309.
\end{document}